\documentstyle[amssymb,amsfonts]{amsart}

\renewcommand{\P}{\mbox{{\bf P}}}

\def\be#1{\begin{equation} \label{#1}}
\def\bi{\begin{itemize}}
\def\bs{\begin{split}}
\def\es{\end{split}}
\def\ba{\begin{align}}
\def\bas{\begin{align*}}
\def\ea{\end{align}}
\def\eas{\end{align*}}
\def\R{{\hbox{\bf R}}}

\def\diag{{\hbox{diag}}}

\def\R{{\hbox{\bf R}}}

\def\div{{\hbox{div}}}
\def\curl{{\hbox{curl}}}
\def\P{{\hbox{\bf P}}}

\def\g{{\frak g}}
\def\eps{\varepsilon}
\def\grad{\nabla}

\def\df{{\rm df}}
\def\cf{{\rm cf}}

\newenvironment{proof}{\noindent {\bf Proof} }{\endprf\par}
\def \endprf{\hfill  {\vrule height6pt width6pt depth0pt}\medskip}
\def\emph#1{{\it #1}}
\def\textbf#1{{\bf #1}}

\theoremstyle{plain}
  \newtheorem{theorem}[subsection]{Theorem}
  
  \newtheorem{proposition}[subsection]{Proposition}
  \newtheorem{lemma}[subsection]{Lemma}

\theoremstyle{remark}

\theoremstyle{definition}

\include{psfig}

\begin{document}

\title[Local well-posedness of Yang-Mills]{Local well-posedness of the Yang-Mills equation in the Temporal Gauge below the energy norm}
\author{Terence Tao}
\address{Department of Mathematics, UCLA, Los Angeles CA 90095-1555}
\email{ tao@@math.ucla.edu}

\vspace{-0.3in}
\begin{abstract}
We show that the Yang-Mills equation in three dimensions in the Temporal gauge is locally well-posed in $H^s$ for $s > 3/4$ if the $H^s$ norm is sufficiently small.  The temporal gauge is slightly less convenient technically than the more popular Coulomb gauge, but has the advantage of uniqueness even for large initial data, and does not require solving a nonlinear elliptic problem.  To handle the temporal gauge correctly we project the connection into curl-free and divergence-free components, and develop some new bilinear estimates of $X^{s,b}$ type which can handle integration in the time direction.
\end{abstract}

\maketitle

\section{Introduction}

This paper is concerned with the low regularity local existence theory for the Cauchy problem for the Yang-Mills equation.  We give only a brief description of this Cauchy problem here; for more detail, see e.g. \cite{kl-mac:yang}.

Let $\g$ be a finite-dimensional Lie algebra, and let $A_\alpha: \R^{3+1} \to \g$ be a $\g$-valued connection on Minkowski space-time, where $\alpha$ ranges over 0,1,2,3.  We use the usual summation conventions on $\alpha$, and raise and lower indices with respect to the Minkowski metric $\eta^{\alpha \beta} := \diag(-1,1,1,1)$.  We define the curvature tensor $F_{\alpha\beta}: \R^{3+1} \to \g$ by
$$ F_{\alpha\beta} := \partial_\alpha A_\beta - \partial_\beta A_\alpha + [A_\alpha,A_\beta]$$
where $[,]$ denotes the Lie bracket on $\g$.  We say that $A_\alpha$ satisfies the \emph{Yang-Mills equation} if $F^{\alpha\beta}$ is divergence-free with respect to the covariant derivative, or more precisely that
$$ \partial_\alpha F^{\alpha\beta} + [A_\alpha, F^{\alpha\beta}] = 0.$$
We can expand this as
$$ \Box A^\beta - \partial^\beta (\partial_\alpha A^\alpha) + \partial_\alpha [A^\alpha, A^\beta]
+ [A_\alpha, \partial^\alpha A^\beta] - [A_\alpha, \partial^\beta A^\alpha]
+ [A_\alpha, [A^\alpha, A^\beta]] = 0$$
where $\Box := \partial_\alpha \partial^\alpha = -\partial_t^2 + \Delta$ is the d'Lambertian.
The Cauchy problem for this equation is not well-posed because of gauge invariance.  However, if one fixes the connection to lie in the \emph{Temporal gauge} $A_0 = 0$, the Yang-Mills equations become essentially hyperbolic, and simplify to 
\be{ym-t}
\partial_t (\div A) + [A_i, \partial_t A_i] = 0
\end{equation}
and
\be{ym-j}
\Box A_j - \partial_j (\div A) + \partial_i [A_i,A_j] + [A_i, \partial_i A_j] - [A_i, \partial_j A_i] + [A_i, [A_i, A_j]] = 0.
\end{equation}
As usual, Roman indices $i,j$ will range over 1,2,3 while Greek indices $\alpha$, $\beta$ range over 0,1,2,3.  Henceforth we shall focus exclusively on the temporal gauge equations \eqref{ym-t}, \eqref{ym-j}, and will no longer consider the situation of more general gauges.

The initial data $A(0), \partial_t A(0)$ to \eqref{ym-t}, \eqref{ym-j} must of course satisfy the compatibility condition
\be{compat}
\div \partial_t A(0) + [A_i(0), \partial_t A_i(0)] = 0.
\end{equation}

Note that one can easily define the notion of a weak solution to these equations \eqref{ym-t}, \eqref{ym-t} as long as $(A, A_t) \in H^s(\R^3) \times H^{s-1}(\R^3)$ for some $s > 1/2$; here $H^s(\R^3)$ denotes the usual Sobolev space $H^s := \{ f: (1-\Delta)^{s/2} f \in L^2 \}$.

The purpose of this paper is to prove 

\begin{theorem}\label{main}  For all $s > 3/4$, the equations \eqref{ym-t}, \eqref{ym-j} are locally well-posed for data in $H^s \times H^{s-1}$ satisfying \eqref{compat} for times $-1 \leq t \leq 1$, if the $H^s \times H^{s-1}$ norm of the data is sufficiently small.  
\end{theorem}

By ``locally well-posed'' we mean that there is a Banach space $X \subseteq C^0_t H^s_x \cap C^1_t H^{s-1}_x$ on the spacetime slab $[-1,1] \times \R^3$ such that for all sufficiently small data $(A, A_t)$ in $H^s \times H^{s-1}$ obeying \eqref{compat}, there is a unique (weak) solution to \eqref{ym-t}, \eqref{ym-j} which lies in the space $X$, and that the map from data to solution is continuous from $H^s \times H^{s-1}$ to $X$.

The analogue of Theorem \ref{main} for the slightly simpler Maxwell-Klein-Gordon equations in the Coulomb gauge was proven in \cite{cuccagna}.  One could also consider the well-posedness of the Yang-Mills equation in the Coulomb gauge, but the resulting equations include a nonlinear elliptic equation which is not always globally solvable, and one must localize the gauge condition somehow.  See \cite{kl-mac:yang} for further discussion.  

The temporal gauge is slightly more difficult technically to handle than the Coulomb gauge, although the two gauges are still quite closely related (for instance, they are more similar to each other than they are to other gauges such as the Lorenz and Cronstrom gauges).  We will emulate the Coulomb gauge analysis (see e.g. \cite{kl-mac:yang}) by splitting the connection $A$ into divergence-free $A^{\df}$ and curl-free $A^{\cf}$ components.  As in the Coulomb gauge, the divergence-free component $A^{\df}$ is the main component of $A$, and has the most interesting dynamics, evolving via a non-linear wave equation \eqref{curl-eq}.  The curl-free component $A^{\cf}$ can be recovered from the divergence-free component by a time integration (see \eqref{div-eq}; contrast with the Coulomb gauge case, where one needs to solve a nonlinear elliptic equation to recover $A_0$ from $A$).  Heuristically, this integration in time should recover one order of smoothness in space; the main technical difficulty is to make this heuristic rigorous and thus allow one to iterate away the nonlinearity.

Theorem \ref{main} could probably be extended in several directions, some of which we discuss below.  However, our purpose here is not so much to obtain a sharp result, but to illustrate that the techniques developed to handle gauge field theories in the Coulomb gauge largely carry over to the Temporal gauge setting, and so allow for an easier treatment of the Yang-Mills equations.  In order to control various time integrals one needs to obtain certain $L^r_x L^q_t$ and $H^s_x H^b_t$ estimates on solutions to the wave equation, which may be of independent interest.

It is likely that Theorem \ref{main} extends to higher dimensions $d > 3$, with the condition on $s$ replaced by $s > \frac{d}{2} - \frac{3}{4}$.  This would be an improvement of $1/4$ a derivative over what can be obtained by Strichartz estimates, although work on simplified models of these equations suggests that one should be able to get within epsilon of the critical regularity (i.e. $s > \frac{d}{2} - 1$) by using more sophisticated function spaces\footnote{Recently, Selberg \cite{selberg:mkg} has been able to achieve this result for $s > \frac{d}{2}-1$ for the closely related Maxwell-Klein-Gordon equation when $d \geq 4$; even more recently, Machedon and Sterbenz \cite{machedon:mkg} have extended this to $d=3$.}, see \cite{kl-mac:null3}, \cite{kl-tar:yang-mills}, \cite{tat:5+1}.  However, the regularity $3/4$ appears to be the best one can do in three dimensions by $X^{s,b}$ type spaces alone.

Since the regularity is sub-critical, the small data assumption should be easily removed by shrinking the time interval.  Unfortunately, the usual technique of exploiting short time intervals by manipulating the $b$ index of the $X^{s,b}$ space (see e.g. \cite{selberg:thesis}) runs into difficulty because there are too many time derivatives on the right-hand side in one of the equations \eqref{div-eq} in the Yang-Mills equation in the Temporal gauge\footnote{In particular, there is not enough room to concede $b$ indices in \eqref{ddc}, \eqref{ddd}.  More precisely, the implicit epsilons in the $b$ exponents $-1/2+$ and $1/2+$ must match, because each derivative in time must reduce the $b$ regularity by at least 1 (since $||\tau|-|\xi||$ is comparable to $|\tau|$ when $|\tau|$ is large).}.  Scaling arguments do not work either, because the $L^2$ component of the $H^s$ norm is super-critical.  This issue appears to be surprisingly delicate.  It appears that one needs to exploit the Lie bracket structure in \eqref{div-eq} in a non-trivial manner to prevent the curl-free portion of the field $A$ from blowing up instantaneously for large data.  One might also need to use the curl-free part of \eqref{ym-j}.

For $s \geq 1$, local and global well-posedness (in either the Temporal or Coulomb gauges) was established in \cite{kl-mac:yang}, with the smooth case achieved earlier in \cite{eardley-moncrief:yang-mills}.  It is plausible that one can adapt the techniques in \cite{keel:mkg} to push the global well-posedness result for Yang-Mills in the Temporal gauge down to $s > 7/8$; however, there is an obstruction because the Hamiltonian does not control the entire $H^1$ norm of the energy in the Temporal gauge\footnote{Using the notation of the sequel, the Hamiltonian is roughly of the form $H(A) \approx \|\nabla A^{\df}\|_2 + \| A^{\cf}_t \|_2 + \|A^{\df}_t\|_2$, and so one loses control of $\| \nabla A^{\cf} \|_2$.  It is possible that one might be able to recover this control by integrating $A^{\cf}_t$ in time (cf. \cite{kl-mac:yang}), but this seems to introduce an additional factor of $T$ for long times $0 \leq t \leq T$, which may worsen the numerology when trying to go below the energy norm.}.

The author thanks Mark Keel and Sergiu Klainerman for helpful discussions, and the anonymous referee for careful reading and many cogent suggestions.  The author also thanks Mark Keel for pointing out a gap in an earlier version of the manuscript, and for suggesting the correct way to fix this gap, and James Grant for additional corrections.  The author is a Clay Prize Fellow supported by grants from the Packard and Sloan foundations.

\section{Notation}\label{notation-sec}

Throughout this paper $s > 3/4$ will be fixed.

We use $A \lesssim B$ to denote the statement that $A \leq CB$ for some constant $C$ depending only on $s$ and the Lie algebra $\g$, and $A \sim B$ to denote the statement $A \lesssim B \lesssim A$.

If $a$ is a number, we use $a+$ to denote a number of the form $a + \eps$ for some $0 < \eps \ll s - 3/4$.  Similarly define $a-$.  We define $\infty-$ to be any sufficiently large number ($\frac{100}{s-3/4}$ will do).

As is usual in the study of gauge field theories (see e.g. \cite{kl-mac:yang}. \cite{kl-mac:null3}) we shall use the projection $\P := \Delta^{-1} (\curl \curl)$ to divergence-free fields, and the companion projection $1 - \P = \Delta^{-1} (\grad \div)$ to curl-free fields.

We use $L^q_t L^r_x$ to denote the space given by the norm
$$ \| u \|_{L^q_t L^r_x} := (\int \| u(t) \|_{L^r_x}^q\ dt)^{1/q}$$
and similarly define $L^r_x L^q_t$, etc. with the obvious modifications when $q = \infty$, or when Lebesgue spaces are replaced by Sobolev spaces $H^s$, etc.

It is by now standard that non-linear wave equations should be studied using the $X^{s,b}$ spaces.  These spaces first appear in \cite{rauch.reed} (see also \cite{beals:xsb}) and were applied to local existence theory by Bourgain, Klainerman and Machedon, and others.  See e.g. \cite{ginibre:survey} for a discussion.  We shall use the notation in \cite{tao:xsb}, and define the wave equation spaces $X^{s,b}_{|\tau| = |\xi|}$ on $\R^3 \times \R$ via the norm
$$ \| u \|_{X^{s,b}_{|\tau| = |\xi|}} :=
\| \langle \xi \rangle^s \langle |\tau| - |\xi| \rangle^b \hat u \|_{L^2_{\xi,\tau}}$$
where the space-time Fourier transform $\hat u$ is defined by
$$ \hat u(\xi,\tau) := \int\int e^{-2\pi i (x \cdot \xi + t \tau)}
u(x,t)\ dx dt$$
and $\langle x \rangle := (1 + |x|^2)^{1/2}$.

We also define the product Sobolev norms $X^{s,b}_{\tau = 0} = H^s_x H^b_t$ by
$$ \| u \|_{X^{s,b}_{\tau = 0}} :=
\| \langle \xi \rangle^s \langle \tau \rangle^b \hat u \|_{L^2_{\xi,\tau}}.$$

In frequency space, the $X^{s,b}_{|\tau| = |\xi|}$ norm is localized to the light cone $|\tau| = |\xi|$ and is thus well adapted for measuring solutions to wave equations $\Box u = F$, while the $X^{s,b}_{\tau = 0}$ is localized to the hyperplane $\tau = 0$  and is thus adapted to measuring solutions to equations such as $\partial_t u = F$.

We use the multiplier norms defined in \cite{tao:xsb}.  Specifically, if $m(\xi_1,\tau_1, \xi_2,\tau_2,\xi_3,\tau_3)$ is a function on the space
$$ \{ (\xi_1,\tau_1,\xi_2,\tau_2,\xi_3,\tau_3) \in (\R^3 \times \R)^3: \xi_1 + \xi_2 + \xi_3 = \tau_1 + \tau_2 + \tau_3 = 0 \}$$
then we define $\|m\|_{[3;\R^3 \times \R]}$ to be the best constant in the inequality
$$ |\int m(\xi_1,\tau_1,\xi_2,\tau_2,\xi_3,\tau_3) \prod_{i=1}^3 \hat u_i(\xi_i,\tau_i) \delta(\xi_1 + \xi_2 + \xi_3) \delta(\tau_1 + \tau_2 + \tau_3)|
\leq \| m \|_{[3;\R^3 \times \R]} \prod_{i=1}^3 \| u_i\|_{L^2_{x,t}}$$
for all $u_1$, $u_2$, $u_3$ on $\R^3 \times \R$.  Similarly define $\| m\|_{[3;\R^3]}$ for purely spatial multipliers.  We shall use several estimates on these multipliers from \cite{tao:xsb} in the sequel.

\section{Preliminary reductions}\label{reductions-sec}

In this section we simplify the equations \eqref{ym-t}, \eqref{ym-j} to a schematic form and reduce matters to proving some bilinear and trilinear $X^{s,b}$ estimates.

The equation \eqref{ym-j} is not yet a non-linear wave equation because of the presence of the gradient $\partial_j (\div A)$.  To eliminate this we shall use the projections $\P$, $(1-\P)$.  More precisely, we split
$$A = A^{\cf} + A^{\df},$$ 
where $A^{\cf} := (1-\P)A$ is the curl-free part of $A$, and $A^{\df} := \P A$ is the divergence-free part of $A$.  The equation \eqref{ym-t} then becomes
\be{div-eq}
\partial_t A^{\cf} = -\Delta^{-1} \grad [A_i, \partial_t A_i]
\end{equation}
while if one applies $\P$ to \eqref{ym-j} one obtains
\be{curl-eq}
\Box A^{\df} = - 2\P[A_i, \partial_i A] + \P[A_i, \grad A_i]
- \P [A_i, [A_i, A]].
\end{equation}

Roughly speaking, the smoothing operator $\Delta^{-1} \nabla$ in \eqref{div-eq} ensures that $A^{\cf}$ has better spatial regularity properties than $A^{\df}$.  In our applications we shall only exploit the fact that $A^{\cf}$ has at least $1/4$ more spatial regularity, but this is not optimal (see \cite{cuccagna} for the analogous situation in Maxwell-Klein-Gordon). 

We shall show that the system \eqref{div-eq}, \eqref{curl-eq}, is locally well-posed for initial data satisfying $\div A^{\df}(0) = \div A^{\df}_t(0) = 0$, $\curl A^{\cf}(0) = 0$, and
\be{acf}
 \| A^{\df}(0)\|_{H^s} + \| A^{\df}_t(0) \|_{H^{s-1}} + \|A^{\cf}(0) \|_{H^s} < \eps
 \end{equation}
for $\eps$ sufficiently small.  Note that one does not need to specify $A^{\cf}_t$ thanks to \eqref{compat} and the fact that \eqref{div-eq} is first-order in time. Note that solutions to the above problem automatically satisfy $\div A^{\df} = \curl A^{\cf} = 0$ for all times $t$ for which the solution exists.  Also, by approximating the data by smooth data and using the global existence and uniqueness of smooth solutions to Yang-Mills in the temporal gauge (\cite{eardley-moncrief:yang-mills}; see also \cite{kl-mac:yang}) we see that the field $A$ constructed in this manner solves \eqref{ym-t}, \eqref{ym-j}, and Theorem \ref{main} follows.

The equation \eqref{acf} only places $A^{\cf}(0)$ in $H^s$.  However, it is possible to take advantage of the gauge symmetry
\begin{equation}\label{gauge-eq} A_i \mapsto U A_i U^{-1} - (\partial_i U) U^{-1}
\end{equation}
of the Yang-Mills system \eqref{ym-t}, \eqref{ym-j} in the temporal gauge to eliminate this quantity\footnote{We thank Mark Keel for this suggestion.}, where $U: \R^3 \to G$ is any field taking values in the (possibly non-compact) Lie group $G$ associated to ${\mathfrak g}$ which does not depend on time (the lack of time dependence is essential in order to preserve the temporal gauge condition $A_0=0$).  Indeed, suppose that 
\begin{equation}\label{a0b}
\|A(0)\|_{H^s} \leq \eps; \quad \|A^{\cf}(0)\|_{H^s} \leq \delta
\end{equation}
for some small $0 < \delta \leq \eps \leq 1$.  From Hodge theory, one can write $A^{\cf}(0) = \div V$ for some potential $V: \R^3 \to \g$ with $\|V\|_{X} \lesssim \delta$, where $\|V\|_X := \| \nabla V \|_{H^s}$ is a partially homogeneous variant of the $H^{s+1}$ norm.  If one then applies the gauge transform \eqref{gauge-eq} with $U := \exp(V)$, then one obtains a new initial data $\tilde A(0)$ of the form 
\begin{equation}\label{tila-eq} \tilde A(0) = A^{\df}(0) + (\exp(V) A^{\df}(0) \exp(-V) - A^{\df}(0)) + (\exp(V) \div V - \div \exp(V)) \exp(-V).
\end{equation}
Using the product estimates
\begin{equation}\label{prod-eq} 
\|f g\|_{X} \lesssim \|f\|_{X} \|g\|_{X}; \quad
\|f g\|_{H^s} \lesssim \|f\|_{X} \|g\|_{H^{s}}; \quad \|f g\|_{H^{s-1}} \lesssim \|f\|_{X} \|g\|_{H^{s-1}}  
\end{equation}
which can be established by standard Littlewood-Paley techniques (the key observation being that $X$ controls the $L^\infty$ norm), and Taylor expansion of $\exp(V), \exp(-V)$, one can verify that the second and third terms on the right-hand side of \eqref{tila-eq} have an $H^s$ norm of $O(\delta \eps)$ if $\eps$ is small enough.  Thus the gauged field $\tilde A(0)$ obeys similar estimates to \eqref{a0b} but with $\eps$ and $\delta$ replaced by $\eps+O(\eps\delta)$ and $O(\eps \delta)$ respectively; also, the gauge transform $U$ used differs from the identity by $O(\delta)$ in $X$ norm.  Iterating this procedure indefinitely starting from data obeying \eqref{acf} and taking limits, we obtain a gauge transform \eqref{gauge-eq} of the initial data for some $U$ with $\| U - 1 \|_{X} \lesssim \eps$, such that the new initial data $\tilde A(0)$ is in the Coulomb gauge\footnote{See also \cite{uhlenbeck} for another construction of the Coulomb gauge which requires less regularity on the initial data, and specifically that $\|A(0)\|_{L^3}$ is sufficiently small.} $\div \tilde A(0) = 0$, and for which \eqref{acf} still holds with $\eps$ replaced by $O(\eps)$.  From the iterative nature of this construction one can also easily verify that the gauge map $A \mapsto \tilde A$ is Lipschitz in the topology given by \eqref{acf}.  Also, from \eqref{prod-eq} we see that applying the gauge transform and its inverse will be continuous in the $C^0_t H^s_x \times C^1_t H^{s-1}_x$ topology.  Because of this, we can now (after redefining $\eps$ slightly) reduce without loss of generality to the case $A^{\cf}(0) = 0$ when we initially in the Coulomb gauge.  (Of course, we do not expect to remain in the Coulomb gauge for non-zero times, as one can already see from \eqref{div-eq}.)  

It thus remains to establish the local well-posedness of the Cauchy problem \eqref{curl-eq}, \eqref{div-eq} with the additional assumption $A^{\cf}(0)=0$.  The equation \eqref{div-eq} we shall treat schematically \footnote{By ``schematically'' we mean that we will ignore such algebraic structures as the Lie bracket, thus for instance $A \partial_t A$ is just shorthand for a tensor of the form $C^i_{jk} A^j \partial_t A^k$ for some constant coefficient tensor $C^i_{jk}$.  We also ignore any Riesz transforms $\Delta^{-1} \nabla^2$, since they are bounded on every space under consideration and so are harmless.  The notation $\nabla^{-1}$ refers to any differential operator of order -1.  As remarked in the introduction, it may be that this algebraic structure needs to be exploited further to obtain large data well-posedness, or to go below 3/4 (cf. \cite{machedon:mkg}).} as
\be{div-eq-simplified}
\partial_t A^{\cf} = \nabla^{-1} (A \partial_t A).
\end{equation}
We could also treat \eqref{curl-eq} schematically as
$$ 
\Box A^{\df} = A \nabla A + A^3$$
but this is a bit too crude, as generic expressions of the form $A^{\df} \nabla A^{\df}$ are too badly behaved (cf. \cite{lindblad:sharpduke}).  To get around this difficulty, we exploit some cancellation in \eqref{curl-eq}.  Specifically, we shall isolate one particular component of \eqref{curl-eq}, namely the \emph{null form} 
\be{null-def}
N(A_1,A_2) := -2\P[(\P A_1)_i, \partial_i A_2] + \P[(A_1)_i, \nabla (A_2)_i].
\end{equation}
The null form is so named because $N(A_1,A_2)$ vanishes whenever $A_1$, $A_2$ are parallel plane waves.  It has the schematic form $N(A_1, A_2) = A_1 \nabla A_2$ but exhibits some additional cancellation.

If we split $A = A^{\df} + A^{\cf}$ in  \eqref{curl-eq}, and use \eqref{null-def} to handle the self-interaction of $A^{\df}$, we may write \eqref{curl-eq} in the schematic form 
\be{curl-eq-simplified}
\Box A^{\df} = N(A^{\df},A^{\df}) + A^{\df} \nabla A^{\cf} + A^{\cf} \nabla A^{\df} + A^{\cf} \nabla A^{\cf} + A^3.
\end{equation}
Note we only need the null structure on the self-interaction of $A^{\df}$; the terms involving $A^{\cf}$ will be better behaved as $A^{\cf}$ will turn out to be 1/4 of a derivative smoother than $A^{\df}$.

To prove the local well-posedness of the system \eqref{curl-eq}, \eqref{div-eq} in $H^s$, we shall iterate in the norm
\be{norm}
\|A\|_X := \| A^{\df} \|_{X^{s,3/4+}_{|\tau| = |\xi|}} +
\| A^{\cf} \|_{X^{s+1/4, 1/2+}_{\tau = 0}};
\end{equation}
this is probably not the only norm for which iteration is possible, but it will suffice for our argument.

From the standard manipulations involving $X^{s,b}$ spaces restricted to fixed time intervals (see e.g. \cite{ginibre:survey}, \cite{selberg:thesis}) it  suffices to estimate $\Box A^{\df}$ in $X^{s-1,-1/4+}_{|\tau| = |\xi|}$ and $\partial_t A^{\cf}$ in $X^{s+1/4, -1/2+}_{\tau = 0}$.  More precisely, by \eqref{curl-eq}, \eqref{div-eq} it suffices to prove the seven estimates
\begin{align}
\| N(A_1,A_2) \|_{X^{s-1,-1/4+}_{|\tau| = |\xi|}}
&\lesssim
\|A_1 \|_{X^{s,3/4+}_{|\tau| = |\xi|}}
\|A_2 \|_{X^{s,3/4+}_{|\tau| = |\xi|}}
\label{null-est}\\
\| A_1 \nabla A_2 \|_{X^{s-1,-1/4+}_{|\tau| = |\xi|}}
+\quad &\nonumber\\
 \| A_2 \nabla A_1 \|_{X^{s-1,-1/4+}_{|\tau| = |\xi|}}
&\lesssim
\|A_1 \|_{X^{s,3/4+}_{|\tau| = |\xi|}}
\|A_2 \|_{X^{s+1/4,1/2+}_{\tau = 0}}
\label{cnd}\\
\| A_1 \nabla A_2 \|_{X^{s-1,-1/4+}_{|\tau| = |\xi|}}
&\lesssim
\|A_1 \|_{X^{s+1/4,1/2+}_{\tau = 0}}
\|A_2 \|_{X^{s+1/4,1/2+}_{\tau = 0}}
\label{dnd}\\
\| A_1 A_2 A_3 \|_{ X^{s-1,-1/4+}_{|\tau| = |\xi|} }
&\lesssim \prod_{i=1}^3 \min( \| A_i \|_{ X^{s,3/4+}_{|\tau| = |\xi|} }, \|A_i \|_{X^{s+1/4,1/2+}_{\tau = 0}} )
\label{cubic}\\
\| \nabla^{-1} (A_1 \partial_t A_2) \|_{X^{s+1/4,-1/2+}_{\tau = 0}}
&\lesssim
\|A_1 \|_{X^{s,3/4+}_{|\tau| = |\xi|}}
\|A_2 \|_{X^{s,3/4+}_{|\tau| = |\xi|}}
\label{cdc}\\
\| \nabla^{-1} (A_1 \partial_t A_2) \|_{X^{s+1/4,-1/2+}_{\tau = 0}}+ \quad&\nonumber\\
\| \nabla^{-1} (A_2 \partial_t A_1) \|_{X^{s+1/4,-1/2+}_{\tau = 0}}
&\lesssim
\|A_1 \|_{X^{s,3/4+}_{|\tau| = |\xi|}}
\|A_2 \|_{X^{s+1/4,1/2+}_{\tau = 0}}
\label{ddc}\\
\| \nabla^{-1} (A_1 \partial_t A_2) \|_{X^{s+1/4,-1/2+}_{\tau = 0}}
&\lesssim
\|A_1 \|_{X^{s+1/4,1/2+}_{\tau = 0}}
\|A_2 \|_{X^{s+1/4,1/2+}_{\tau = 0}}
\label{ddd}
\end{align}
for all fields $A_1, A_2, A_3$ on $\R^{3+1}$ (not necessarily divergence-free or curl-free).  Indeed, if the estimates \eqref{null-est}-\eqref{ddd} held, then we would have the a priori estimate
$$
\| \Box A^{\df} \|_{X^{s-1,-1/4+}_{|\tau| = |\xi|}} +
\| \partial_t A^{\cf} \|_{X^{s+1/4, -1/2+}_{\tau = 0}}
\lesssim \|A\|_X^2 + \|A\|_X^3$$
for solutions to \eqref{curl-eq}, \eqref{div-eq}, where $X$ was the norm in \eqref{norm}.  From the standard energy estimates\footnote{The standard energy estimates actually give a little extra regularity on $A^{\df}$, allowing one to obtain control on $\partial_t A^{\df}$.  Unfortunately we do not have a similar amount of surplus time regularity for $A^{\cf}$, which is why we were unable to extend this result to large data.} for $X^{s,b}$ spaces (see \cite{ginibre:survey}, \cite{selberg:thesis}), and the smallness propreties \eqref{acf} of the initial data, and the initial Coulomb gauge condition $A^\cf(0) = 0$, we thus have
$$ \|A\|_X \lesssim \eps + \|A\|_X^2 + \|A\|_X^3,$$
which will give an a priori bound on $\|A\|_X$ by continuity arguments if $\eps$ is sufficiently small.  One can then set up a standard Picard iteration scheme (see e.g. \cite{ginibre:survey},\cite{selberg:thesis}) in the Banach space $X$ for the Cauchy problem \eqref{curl-eq}, \eqref{div-eq} and adapt the above argument to differences of solutions to obtain local well-posedness (indeed one even obtains analytic dependence of the solution on the initial data this way).

From \eqref{null-def} it is well-known (see \cite{kl-mac:yang}, \cite{kl-tar:yang-mills}) that $N$ is of the schematic form
$$
N(A_1,A_2) := Q(\nabla^{-1} A_1, A_2) + \nabla^{-1} Q(A_1,A_2)
$$
where $\nabla^{-1}$ is some Fourier multiplier of order $-1$, and $Q$ is some finite linear combination of the null forms
$$ Q_{ij}(A_1, A_2) := \partial_{x_i} A_1 \partial_{x_j} A_2 - \partial_{x_j} A_1 \partial_{x_i} A_2.$$
The claim \eqref{null-est} then follows immediately from \cite{tao:xsb}, Proposition 9.2 (or \cite{cuccagna}, Lemma 4; see also \cite{keel:mkg}).

We remark that the main difference between the full Yang-Mills equation and the simplified models studied in e.g. \cite{kl-tar:yang-mills}, for the purposes of local existence theory, is that the treatment of the model only requires estimates similar to \eqref{null-est}, and not the additional estimates \eqref{cnd}-\eqref{ddd} arising from the gauge\footnote{In \cite{kl-tar:yang-mills} the regularity is so close to critical that one cannot rely purely on $X^{s,b}$ spaces alone, and must iterate in more complicated spaces which have a physical space component in addition to a frequency space component even for the model equation.  Nevertheless, recent work in \cite{selberg:mkg}, \cite{machedon:mkg} has been able to use these more complicated norms for genuine gauge equations, and not just the model equations.}.  Note that the new estimates \eqref{cnd}-\eqref{ddd} which need to be proven do not require any sort of null structure.

It of course remains to prove \eqref{cnd}-\eqref{ddd}.  These are hybrid estimates combining the wave $X^{s,b}_{|\tau|=|\xi|}$ spaces with the product Sobolev spaces $X^{s,b}_{\tau=0}$, and so are not covered by systematic tables of estimates such as those in \cite{damiano:null} or \cite{tao:xsb}.  These hybrid estimates can be proven solely by a large number of applications of the Cauchy-Schwarz inequality, however we have elected to use a mix of techniques, combining the machinery of \cite{tao:xsb} with some Strichartz estimates which may be of independent interest.  We shall also use some algebraic ``denominator games'' to redistribute various Fourier weights; these can be viewed as variants of the fractional Leibnitz rule.  (An equivalent approach would have been to decompose the multiplier into various regions such as $|\xi_1| \sim |\xi_2| \gtrsim |\xi_3|$ and then simplify the weights on each region individually).

For a first reading of the proofs of \eqref{cnd}-\eqref{ddd} we recommend setting $s=3/4$ and ignoring epsilons.

In the estimates \eqref{cdc}-\eqref{ddd} the expressions $\| \nabla^{-1} (\ldots) \|_{ X^{s+1/4,-1/2+}_{\tau = 0} }$ can automatically be replaced by the slightly smaller quantities $\| \ldots \|_{X^{s-3/4,-1/2+}_{\tau = 0}}$ by general separation of scale arguments (see\footnote{More precisely, one can write all the estimates \eqref{cdc}-\eqref{ddd} as bounds on multiplier norms $\| m \|_{[3; \R^3 \times \R]}$ or $\| m \|_{[4; \R^3 \times \R]}$ for various multipliers $m$, which in the $\xi$ variables have singularities which blow up like $1/|\xi_j|$ at worst.  Corollary 8.2 in \cite{tao:xsb} then allows us to replace those singularities with $1/\langle \xi_j \rangle$, which is equivalent to the replacement mentioned above.  Intuitively, the explanation for this is that low frequencies $|\xi| \ll 1$ do not play a significant role in the short-time existence theory because the uncertainty principle does not give them enough time to properly form for times $|t| \leq 1$; note that in wave equations time and space have a similar scaling.} \cite{tao:xsb}, Corollary 8.2).  We shall implicitly assume this replacement in the sequel. 

\section{Strichartz estimates}\label{strichartz-sec}

In this section we list some estimates of Strichartz type which will be useful in proving \eqref{cnd}-\eqref{ddd}.

We have the well-known energy estimate
\be{energy-est}
\| u \|_{L^\infty_t H^s_x} \lesssim \| u \|_{X^{s,1/2+}_{|\tau| = |\xi|}}.
\end{equation}
We also have the Strichartz estimate
\be{strich}
\| u \|_{L^4_{t,x}} \lesssim \| u \|_{X^{1/2,1/2+}_{|\tau| = |\xi|}};
\end{equation}
however for our purposes we shall need the bilinear improvement
\be{bilinear-strichartz}
\| uv \|_{L^2_{t,x}} \lesssim \| u \|_{X^{1/2+\delta,1/2+}_{|\tau| = |\xi|}}
\| v \|_{X^{1/2-\delta,1/2+}_{|\tau| = |\xi|}}
\end{equation}
for all $|\delta| \leq 1/2$, see \cite{damiano:null}, \cite{kl-mac:null} (or \cite{borg:book} for the analogous estimate for Schr\"odinger).  Actually we shall only need this estimate with $|\delta| < 1/4$.

We shall also use another estimate of Strichartz type (but with $L^r_x L^q_t$ norms rather than $L^q_t L^r_x$), which does not seem to be explicitly in the literature\footnote{A similar estimate however was used by Tataru \cite{kl-mac:null}.}:

\begin{proposition}\label{secret-strichartz}
For any $u$ in $\R^{3+1}$, we have
$$ 
\| u \|_{L^4_x L^2_t} \lesssim \| u \|_{X^{1/4, 1/2+}_{|\tau| = |\xi|}}.$$
\end{proposition}

\begin{proof}
By the usual averaging over time modulations argument\footnote{See e.g. \cite{ginibre:survey}, \cite{selberg:thesis}.  The idea is to use the identity $u(t,x) = \int_\R e^{it\lambda} u_\lambda(t,x)\ d\lambda$, where $\widehat{u_\lambda}(\tau,\xi) := \delta(\pm\tau-|\xi|) \hat u(\tau + \lambda,\xi)$ and $\tilde u$ is supported on the half-plane $\pm \tau \geq 0$.  This expresses a general $X^{1/4,1/2+}$ function as a modulated average of $H^{1/4}$ free solutions.  Since the $e^{it\lambda}$ phase has no impact on the $L^4_x L^2_t$ norm, the reduction then follows from Minkowski's inequality and Cauchy-Schwarz (noting that the weight $\langle \lambda \rangle^{-1/2-}$ is in $L^2_\lambda$).} we may assume that $u$ is a free solution to the wave equation.  By time reversal symmetry we may assume $u$ is a forward solution
$$ u(t) = e^{it\sqrt{-\Delta}} u_0$$
in which case it suffices to show the scale-invariant estimate
$$ \| e^{it\sqrt{-\Delta}} u_0 \|_{L^4_x L^2_t} \lesssim \| u_0 \|_{\dot H^{1/4}}.$$
By the usual Littlewood-Paley arguments (e.g. \cite{sogge:wave}) we may assume that $u_0$ is restricted to a dyadic annulus; by scale invariance we may assume that $\hat u_0$ is supported on the annulus $|\xi| \sim 1$.

By Plancherel in $t$ we have
$$ \| e^{it\sqrt{-\Delta}} u_0 \|_{L^4_x L^2_t} \sim \| \delta(\tau - \sqrt{-\Delta}) u_0 \|_{L^4_x L^2_\tau}.$$
By Plancherel in $x$, polar co-ordinates, and the frequency localization we have
$$
\| u_0 \|_{\dot H^{1/4}} \sim \| u_0 \|_2 \sim \| (\int_{|\xi| = \tau} |\hat u(\xi)|^2\ d\xi)^{1/2} \|_{L^2_\tau}.$$
By the frequency localization we may restrict $\tau$ to the range $\tau \sim 1$.  The claim then follows from the embedding $L^2_\tau L^4_x \subset L^4_x L^2_\tau$ and the Stein-Tomas-Sj\"olin $(L^2,L^4)$ restriction theorem 
$$ \| \delta(\tau - \sqrt{-\Delta}) u \|_{L^4_x} \lesssim (\int_{|\xi| = \tau} |\hat u(\xi)|^2\ d\xi)^{1/2} \hbox{ for } \tau \sim 1$$
for the sphere $S^2$ (see e.g. \cite{stein:large}).
\end{proof}

This Proposition should be compared (using Sobolev embedding) with \eqref{strich} as well as with the estimate
$$ \| u \|_{L^\infty_x L^2_t} \lesssim \| u \|_{X^{1, 1/2+}_{|\tau| = |\xi|}}$$
proven in \cite{kl-mac:yang}.  

\section{The cubic term: Proof of \eqref{cubic}}\label{cubic-sec}

In the remainder of this paper we prove the estimates \eqref{cnd}-\eqref{ddd}.  We shall tackle these estimates in increasing order of difficulty, beginning with the easy cubic estimate \eqref{cubic}.

There is plenty of room in this estimate, and we shall be somewhat generous with regularity.

By interpolating \eqref{energy-est} for $s=0$ with the trivial identity
$X^{0,0}_{|\tau| = |\xi|} = L^2_t L^2_x$ we have
$$
\| u \|_{L^{4-}_t L^2_x} \lesssim \| u \|_{X^{0, 1/4-}_{|\tau| = |\xi|}}.
\lesssim \| u \|_{X^{1-s, 1/4-}_{|\tau| = |\xi|}}
$$
so by duality
$$
\| u \|_{X^{s-1, -1/4+}_{|\tau| = |\xi|}} \lesssim \| u \|_{L^{4/3+}_t L^2_x}.
$$
To show \eqref{cubic} it thus suffices by H\"older to show that
$$
\| u \|_{L^{4+}_t L^6_x} \lesssim \| u \|_{X^{s,3/4+}_{|\tau| = |\xi|}}
$$
and
$$
\| u \|_{L^{4+}_t L^6_x} \lesssim \| u \|_{X^{s+1/4,1/2+}_{\tau = 0}}.
$$
The former follows from \eqref{strich} and Sobolev, since $s > 3/4$.  The latter follows from the Sobolev embeddings $H^{1/2+}_t \subset L^{4+}_t$
and $H^{s+1/4}_x \subset H^1_x \subset L^6_x$.
This shows \eqref{cubic}.

\section{Curl-free interactions: Proof of \eqref{ddd}, \eqref{dnd}}

In this section we prove \eqref{ddd}, \eqref{dnd}, which control the interactions between two curl-free fields.  These estimates are relatively easy because the curl-free component is quite regular; indeed, they will mostly follow from Sobolev embedding and H\"older.  We first show the preliminary estimates

\begin{lemma}\label{trio}
We have
\be{trio-1}
|\int\int uvw\ dx dt| \lesssim
\| u \|_{X^{3/4 - s,1/2-}_{\tau = 0}}
\| v \|_{X^{s+1/4,1/2+}_{\tau = 0}}
\| w \|_{X^{s+1/4,-1/2+}_{\tau = 0}}.
\end{equation}
and
\be{trio-2}
|\int\int uvw\ dx dt| \lesssim
\| u \|_{X^{s-3/4,1/2+}_{\tau = 0}}
\| v \|_{X^{s+1/4,1/2+}_{\tau = 0}}
\| w \|_{X^{5/4-s-,-1/4+}_{\tau = 0}}.
\end{equation}
\end{lemma}

\begin{proof}
By \cite{tao:xsb}, Lemma 3.6 the above estimates split into the spatial estimates
\begin{align*}
|\int\int fgh\ dx| &\lesssim
\| f \|_{H^{3/4-s}_x(\R^3)}
\| g \|_{H^{s+1/4}_x(\R^3)}
\| h \|_{H^{s+1/4}_x(\R^3)}\\
|\int\int fgh\ dx| &\lesssim
\| f \|_{H^{s-3/4}_x(\R^3)}
\| g \|_{H^{s+1/4}_x(\R^3)}
\| h \|_{H^{5/4-s-}_x(\R^3)}
\end{align*}
and the temporal estimates
\begin{align*}
|\int\int FGH\ dx| &\lesssim
\| F \|_{H^{1/2-}_t(\R)}
\| G \|_{H^{1/2+}_t(\R)}
\| H \|_{H^{-1/2+}_t(\R)}\\
|\int\int FGH\ dx| &\lesssim
\| F \|_{H^{1/2+}_t(\R)}
\| G \|_{H^{1/2+}_t(\R)}
\| H \|_{H^{-1/4+}_t(\R)}.
\end{align*}
But these follow from Sobolev multiplication laws (\cite{tao:xsb}, Proposition 3.15, or \cite{selberg:thesis}; alternatively, use fractional Leibnitz, Sobolev embedding and H\"older) and the hypothesis $s > 3/4$.  (In fact, one has some regularity to spare in all of these estimates, except for the first temporal estimate).
\end{proof}

By duality and the remarks at the end of Section \ref{reductions-sec}, 
\eqref{ddd} follows from \eqref{trio-1} and the easily verified embedding
$$
\| w_t \|_{X^{s+1/4,-1/2+}_{\tau = 0}} \lesssim \| w \|_{X^{s+1/4,1/2+}_{\tau = 0}}.
$$

The estimate \eqref{dnd} similarly follows from duality, \eqref{trio-2}, and the easily verified embeddings
$$
\| \nabla u \|_{ X^{s-3/4,1/2+}_{\tau = 0}} \lesssim \| u \|_{ X^{s+1/4,1/2+}_{\tau = 0} }
$$
and
$$
\| w \|_{X^{5/4-s-,-1/4+}_{\tau = 0}} \lesssim \| w \|_{X^{1-s,1/4-}_{|\tau| = |\xi|}}
$$
(the latter following from the crude estimate $\langle \xi \rangle / \langle \tau \rangle \lesssim \langle |\tau| - |\xi| \rangle$).

\section{Hybrid curl-free interactions: Proof of \eqref{ddc}}\label{hybrid-sec}

We now prove \eqref{ddc}, which controls the extent to which a curl-free and div-free field may interact and produce another curl-free field.

We first observe that we may assume that $\hat A_1$ is restricted to the neighbourhood of the light cone
\be{near-light}
||\tau| - |\xi|| \ll |\xi|.
\end{equation}
This is because if $\hat A_1$ vanishes on \eqref{near-light}, then 
we have the embedding
\be{c-embed}
\|A_1 \|_{X^{s+1/4,1/2+}_{\tau = 0}} \lesssim \| A_1 \|_{X^{s,3/4+}_{|\tau| = |\xi|}}
\end{equation}
and so \eqref{ddc} would then follow from \eqref{ddd}, which was previously proven.

By duality it thus suffices to show the estimates
$$
|\int\int u v_t w\ dx dt| + |\int\int u v w_t\ dx dt| 
\lesssim 
\| u \|_{X^{3/4-s,1/2-}_{\tau = 0}}
\| v \|_{X^{s+1/4,1/2+}_{\tau = 0}}
\| w \|_{X^{s,3/4+}_{|\tau| = |\xi|}}
$$
whenever $\hat w$ is supported on \eqref{near-light}.

Using the notation of \cite{tao:xsb}, it suffices to show that
$$
\| \frac{ (|\tau_2| + |\tau_3|) \chi_{||\xi_3|-|\tau_3|| \ll |\xi_3|}}{
\langle \xi_1 \rangle^{3/4-s} \langle \tau_1 \rangle^{1/2-}
\langle \xi_2 \rangle^{s+1/4} \langle \tau_2 \rangle^{1/2+}
\langle \xi_3 \rangle^{s} \langle |\tau_3| - |\xi_3| \rangle^{3/4} } \|_{[3;\R^3 \times \R]} \lesssim 1.$$
Since $||\xi_3|-|\tau_3|| \ll |\xi_3|$, we have $\langle \tau_3 \rangle \sim \langle \xi_3 \rangle$.  Since $\tau_1 + \tau_2 + \tau_3 = 0$, we thus have
$$ |\tau_2| + |\tau_3| 
\lesssim 
\langle \tau_1 \rangle^{1/2-} \langle \tau_2 \rangle^{1/2+} +
\langle \tau_1 \rangle^{1/2-} \langle \xi_3 \rangle^{1/2+}
+ \langle \tau_2 \rangle^{1/2+} \langle \xi_3 \rangle^{1/2-}.$$
Inserting this into the above the Comparison principle (\cite{tao:xsb}, 
Lemma 3.1) and rewriting things in integral form, we reduce to showing the 
estimates
\begin{align}
|\int\int u v w\ dx dt| 
 & \lesssim 
\| u \|_{X^{3/4-s,0}_{\tau = 0}}
\| v \|_{X^{s+1/4,0}_{\tau = 0}}
\| w \|_{X^{s,3/4+}_{|\tau| = |\xi|}}
\label{f1}\\
|\int\int u v w\ dx dt| 
 & \lesssim 
\| u \|_{X^{3/4-s,0}_{\tau = 0}}
\| v \|_{X^{s+1/4,1/2+}_{\tau = 0}}
\| w \|_{X^{s-1/2-,3/4+}_{|\tau| = |\xi|}}
\label{f2}\\
|\int\int u v w\ dx dt| 
 & \lesssim 
\| u \|_{X^{3/4-s,1/2-}_{\tau = 0}}
\| v \|_{X^{s+1/4,0}_{\tau = 0}}
\| w \|_{X^{s-1/2-,3/4+}_{|\tau| = |\xi|}}.
\label{f3}
\end{align}

To show \eqref{f1}, it suffices by \eqref{energy-est} and the hypothesis $s > 3/4$ to show that (conceding some derivatives)
$$
|\int\int u v w\ dx dt| \lesssim \| u \|_{L^2_{t,x}} \| v \|_{L^2_t H^1_x}
\| w \|_{L^\infty_t H^{1/2}_x}.$$
But this follows from Sobolev and H\"older.

To show \eqref{f2}, \eqref{f3}, it suffices by Proposition \ref{secret-strichartz} and the hypothesis $s > 3/4$ to show that (conceding some derivatives)
$$
|\int\int u v w\ dx dt| \lesssim \| u \|_{L^2_{t,x}} \| v \|_{H^{1/2+}_t H^{3/4}_x}
\| w \|_{L^2_t L^4_x}$$
and
$$
|\int\int u v w\ dx dt| \lesssim \| u \|_{H^{1/2+}_t L^2_x} \| v \|_{L^2_t H^{3/4}_x}
\| w \|_{L^2_t L^4_x}.$$
But these estimates follow from Sobolev and H\"older.  This completes the proof of \eqref{ddc}.

\section{Divergence-free interactions: Proof of \eqref{cdc}}\label{cdc-sec}

We now prove \eqref{cdc}, which controls how two div-free fields may interact to cause a curl-free field.  We may assume as before that $\hat A_1$, $\hat A_2$ are supported in the region \eqref{near-light}, since the claim follows from \eqref{c-embed} and the previously proven estimates \eqref{ddc}, \eqref{ddd} otherwise.  By duality we can thus rewrite the estimate as
$$
\| \frac{ |\tau_3| \chi_{||\xi_2|-|\tau_2|| \ll |\xi_2|}
\chi_{||\xi_3|-|\tau_3|| \ll |\xi_3|}}{
\langle \xi_1 \rangle^{3/4-s} \langle \tau_1 \rangle^{1/2-}
\langle \xi_2 \rangle^{s} \langle |\tau_2| - |\xi_2| \rangle^{3/4+}
\langle \xi_3 \rangle^{s} \langle |\tau_3| - |\xi_3| \rangle^{3/4+} } \|_{[3;\R^3 \times \R]} \lesssim 1.$$
Since $||\xi_i| - |\tau_i|| \ll |\xi_i|$ for $i=2,3$ and $\tau_1 + \tau_2 + \tau_3 = 0$, we have the estimate
$$ |\tau_3| \lesssim 
\langle \tau_1 \rangle^{1/2-} \langle \xi_3 \rangle^{1/2+}
+ \langle \xi_2 \rangle^{1/2-} \langle \xi_3 \rangle^{1/2+}.$$
Inserting this into the previous using the Comparison principle (\cite{tao:xsb}, Lemma 3.1) and rewriting things in integral form, we reduce to proving 
\begin{align}
|\int\int u v w\ dx dt| 
 & \lesssim 
\| u \|_{X^{3/4-s,0}_{\tau = 0}}
\| v \|_{X^{s-1/2+,3/4+}_{|\tau| = |\xi|}}
\| w \|_{X^{s,3/4+}_{|\tau| = |\xi|}}
\label{g1}\\
|\int\int u v w\ dx dt| 
 & \lesssim 
\| u \|_{X^{3/4-s,1/2-}_{\tau = 0}}
\| v \|_{X^{s-1/2+,3/4+}_{|\tau| = |\xi|}}
\| w \|_{X^{s-1/2-,3/4+}_{|\tau| = |\xi|}}
\label{g2}
\end{align}

To prove \eqref{g1}, it suffices from the hypothesis $s > 3/4$ and Cauchy-Schwarz to show that
$$ \| vw\|_{L^2_{x,t}} \lesssim
\| v \|_{X^{s-1/2+,3/4+}_{|\tau| = |\xi|}}
\| w \|_{X^{s,3/4+}_{|\tau| = |\xi|}}.$$
But this follows from \eqref{bilinear-strichartz} and the hypothesis $s > 3/4$.

To prove \eqref{g2}, we observe from Proposition \ref{secret-strichartz}, the hypothesis $s > 3/4$ and a little bit of Sobolev embedding that
$$ \| v \|_{L^{2+}_t L^4_x} \lesssim \| v \|_{X^{s-1/2 \pm,3/4+}_{|\tau| = |\xi|}}$$
and similarly from $w$.  The claim then follows from H\"older and the Sobolev embedding
$$ \| u \|_{L^{\infty-}_t L^2_x} \lesssim \| u \|_{X^{3/4-s,1/2-}_{\tau = 0}}.$$

\section{Hybrid divergence-free interactions: Proof of \eqref{cnd}}\label{cnd-sec}

We now prove \eqref{cnd}, which controls how a curl-free and div-free field may interact to cause a div-free field; this is the most difficult of all the estimates.

We can rewrite \eqref{cnd} as
$$
\| \frac{
(|\xi_2| + |\xi_3|)
\langle \xi_1 \rangle^{s-1} \langle |\tau_1| - |\xi_1| \rangle^{-1/4+}
}
{
\langle \xi_2 \rangle^s
\langle |\xi_2| - |\tau_2| \rangle^{3/4+}
\langle \xi_3 \rangle^{s+1/4}
\langle \tau_3 \rangle^{1/2+}
}
\|_{[3; \R^3 \times \R]} \lesssim 1.$$

We may restrict to the region $|\xi_2| \leq |\xi_1|$ since the other region then follows by symmetry and the Comparison principle (\cite{tao:xsb}, Lemma 3.1).  In this case we may estimate $|\xi_2| + |\xi_3|$ by $\langle \xi_1 \rangle$.

We will discard the bounded factor $\langle |\tau_1| - |\xi_1| \rangle^{-1/4+}$.  By two applications of the averaging argument in (\cite{tao:xsb}, Proposition 5.1) it then suffices to show that
$$
\| \frac{
\langle \xi_1 \rangle^s 
\chi_{||\xi_2| - |\tau_2|| \sim 1}
\chi_{|\tau_3| \sim 1}
}
{
\langle \xi_2 \rangle^s
\langle \xi_3 \rangle^{s+1/4}
}
\|_{[3; \R^3 \times \R]} \lesssim 1.$$
Suppose we restrict $\tau_2$ to the region $\tau_2 = T + O(1)$ for some integer $T$, then $\tau_1$ is then restricted to the region $\tau_1 = -T + O(1)$, and $\xi_2$ is restricted to the annulus $|\xi_2| = |T| + O(1)$.  The $\tau_1$ regions are essentially disjoint as $T$ varies along the integers, and similarly for $\tau_2$.  By Schur's test (\cite{tao:xsb}, Lemma 3.11) it thus suffices to show that
$$
\| \frac{
\langle \xi_1 \rangle^s 
\chi_{\tau_1 = -T + O(1)}
\chi_{\tau_2 = T + O(1)}
\chi_{|\xi_2| = |T| + O(1) }
\chi_{|\tau_3| \sim 1}
}
{
\langle \xi_2 \rangle^s
\langle \xi_3 \rangle^{s+1/4}
}
\|_{[3; \R^3 \times \R]} \lesssim 1$$

uniformly in $T$.  The $\tau$ behaviour is now trivial, and we can reduce (e.g. by \cite{tao:xsb} Lemmata 3.6 and 3.14) to the spatial estimate
$$
\| \frac{
\langle \xi_1 \rangle^s 
\chi_{|\xi_2| = |T| + O(1) }
}
{
\langle T \rangle^s
\langle \xi_3 \rangle^{s+1/4}
}
\|_{[3; \R^3]} \lesssim 1.$$
We may of course assume that $T$ is a positive integer.  We may assume that $|\xi_3| \leq |\xi_1|$, since the other case then follows by symmetry and the Comparison principle (\cite{tao:xsb}, Lemma 3.1).

There are only two remaining cases in which the symbol does not vanish: $|\xi_1| \sim |\xi_3| \gtrsim T$ and $|\xi_1| \sim T \gtrsim |\xi_3|$.  In the first case we reduce to
$$
\| \frac{ \chi_{|\xi_2| = T + O(1) }
}
{
T^{s+1/4}
}
\|_{[3; \R^3]} \lesssim 1.$$
On the other hand, from a Cauchy-Schwarz estimate (\cite{tao:xsb}, Lemma 3.14) 
we have
$$
\| \chi_{|\xi_2| = T + O(1)} \|_{[3;\R^3]} \leq |\{ \xi_2: |\xi_2| = T + O(1) \}|^{1/2} \lesssim T,$$
and the claim follows since $s > 3/4$.

Now suppose that $|\xi_1| \sim T \gtrsim |\xi_3|$.  We thus reduce to
$$
\| 
\chi_{|\xi_2| = T + O(1) }
\langle \xi_3 \rangle^{s+1/4}
\|_{[3; \R^3]} \lesssim 1.$$
At this point, we resort to Cauchy-Schwarz (\cite{tao:xsb}, Corollary 3.10), and reduce to showing that
$$ \| \chi_{|\xi| = T + O(1)} * \langle \xi \rangle^{2s+1/2} \|_\infty \lesssim 1.$$
But this is easily verified since $s > 3/4$.

\end{document}